\documentclass[12pt]{article}
\usepackage{graphicx,amsmath,amsfonts,amssymb,amscd,bezier,amstext,makeidx}
\usepackage{fancyheadings}
\usepackage{amsthm}
\usepackage{bm}
\usepackage{color}
\usepackage{fancybox}
\usepackage{fancyheadings}
\usepackage{latexsym}
\usepackage{graphicx}
\usepackage{amsthm,amssymb,amscd}
\usepackage[all]{xy}
\usepackage{longtable}

\begin{document}

\newtheorem{definition}{Definition}[section]
\newtheorem{theorem}{Theorem}[section]
\newtheorem{corollary}{Corollary}[section]
\newtheorem{lemma}{Lemma}[section]
\newtheorem{proposition}{Proposition}[section]
\newtheorem{example}{Example}[section]
\newtheorem{remark}{Remark}[section]

\newfont{\bms}{msbm10 scaled 1200}
\def\matR{\mbox{\bms R}}
\def\matF{\mbox{\bms F}}
\def\matZ{\mbox{\bms Z}}
\def\matN{\mbox{\bms N}}
\def\matQ{\mbox{\bms Q}}
\def\matK{\mbox{\bms K}}
\def\matL{\mbox{\bms L}}
\def\matC{\mbox{\bms C}}
\def\matM{\mbox{\bms M}}
\newcommand{\B}{b\!\!\!b}
\newcommand{\A}{e\!\!\!a}

\title{Codimension Two Determinantal Varieties with Isolated Singularities}

\author{{ Miriam da Silva Pereira }
\thanks{Partially supported by FAPESP grant \# 05/58960-3.}  \\
{\small $$}{\footnotesize Centro de Ci\^{e}ncias Exatas e da Natureza- UFPB}\\
{\small $$}{\footnotesize Cidade Universit\'aria}\\
{\small $$}{\footnotesize e-mail: miriam@mat.ufpb.br - Jo\~{a}o
Pessoa - PB Brazil}
\\{\small $$} \\ {Maria Aparecida Soares Ruas}
\thanks{Partially supported by CNPq grant \# 303774/2008-8 and FAPESP
grant \# 08/54222-6.}  \\
{\small $$}{\footnotesize Instituto de Ci\^{e}ncias Matem\'aticas e de Computa\c c\~ao- USP}\\
{\small $$}{\footnotesize Avenida Trabalhador S\~{a}o-carlense, 400}\\
{\small$$}{\footnotesize e-mail: maasruas@icmc.usp.br - S\~{a}o
Carlos - SP Brazil} }

\maketitle

\begin{abstract}
\noindent We study codimension two determinantal varieties with
isolated singularities. These singularities admit a unique
smoothing, thus we can define their Milnor number as the middle
Betti number of their generic fiber. For surfaces in $\matC^4$, we
obtain a L\^{e}-Greuel formula expressing the Milnor number of the
surface in terms of the second polar multiplicity and the Milnor
number of a generic section. We also relate the Milnor number with
Ebeling and Gusein-Zade index of the $1$- form given by the
differential of a generic linear projection defined on the surface.
To illustrate the results, in the last section we compute the Milnor
number of some normal forms from A. Fr\"uhbis-Kr\"uger and A. Neumer
\cite{FKn} list of simple determinantal surface singularities.
\end{abstract}


\section{Introduction}

The goal of this paper is to study codimension two determinantal
varieties with an isolated singularity. These singularities admit a
unique smoothing, hence  the topological type of their Milnor fiber
is well defined.

Let $X$ be a codimension two determinantal variety with isolated
singularity and $X_t$ its generic fiber. We define the Milnor number
of  $X$ as the middle Betti number of $X_t$. The condition that $X$
has isolated singularity implies that $dim(X)=2,\,3$. Since these
are normal singularities, it follows from a result of Greuel and
Steenbrink ( \cite{greuel}, pg. 540) that $b_1(X_t)=0$, where $b_1$
is the first Betti number. They also prove in \cite{greuel} that for
every complex analytic space with isolated singularity one has
$\pi_i(X_t)=0$, for $i\leq dim(X)-codim(X)$, where $\pi_i(X_t)$ is
the $i$-th homotopy group of $X_t$. Thus, it follows that the
generic fiber of a determinantal variety $X_t$ with isolated
singularity is connected. When $dim(X)=3$, it also follows that
$X_t$ is $1$-connected.

For determinantal surfaces $X$ in $\matC^4$, we use these results
and Morse theory to obtain a L\^{e}-Greuel formula expressing the
Milnor number $\mu(X)$, in terms of the second polar multiplicity
$m_2(X)$ and the Milnor number of a generic section of $X$. This
formula holds for $3$- dimensional determinantal varieties in
$\matC^5$, under the adicional hypothesis that $b_2(X_t)=0$.

We do not know an algebraic formula to compute $\mu(X_t)$. Our
approach in this paper, in order to calculate this invariant is to
further investigate its geometric interpretation. For this we relate
$m_2(X)$, and consequently $\mu(X)$, to the Ebeling and Gusein-Zade
index of the $1$- form $dp$, where $p$ is a generic linear
projection defined on $X$.

We show in the last section how to use the results to compute the
Milnor number of some normal forms from A. Fr\"uhbis-Kr\"uger and A.
Neumer \cite{FKn} list of simple determinantal surface
singularities.

The recent paper \cite{BOT} by Ballesteros, Or\'efice and Tomazella
also discusses the Milnor number of functions on determinantal
varieties.


\section{Determinantal Varieties }

Let $Mat_{(n,p)}(\matC)$ be the set of all $n\times p$ matrices with
complex entries, $\Delta_t\subset Mat_{(n,p)}(\matC)$ the subset
formed by matrices that have rank less than $t$, with $1\leq t\leq
\min(n,p)$. It is possible to show that $\Delta_t$ is an irreducible
singular algebraic variety of codimension $(n-t+1)(p-t+1)$ (see
\cite{Bruns}). Moreover the singular set of $\Delta_t$ is exactly
$\Delta_{t-1}$. The set $\Delta_t$ is called \textit{generic
determinantal variety}.

\begin{definition}
Let $M=(m_{ij}(x))$ be a $n\times p$ matrix whose entries are
complex analytic functions on $U\subset\matC^r$, $0\in U$ and $f$
the function defined by the $t\times t$ minors of $M$. We say that
$X$ is a determinantal variety of codimension $(n-t+1)(p-t+1)$ if
$X$ is defined by the equation $f=0$.
\end{definition}

We can look to a matrix $M=(m_{ij}(x))$ as a map
$M:\matC^r\longrightarrow Mat_{(n,p)}(\matC)$, with $M(0)=0$. Then,
the determinantal variety in $\matC^r$ is the set
$X=M^{-1}(\Delta_t)$, with $1\leq t\leq \min\{n,p\}$. The singular
set of $X$ is given by $M^{-1}(\Delta_{t-1})$. We denote
$X_{reg}=M^{-1}(\Delta_t\backslash \Delta_{t-1})$, the regular part
of $X$. Notice that $X$ has isolated singularity at the origin if
and only if $r\leq(n-t+2)(p-t+2)$.

Let ${\cal{O}}_r$ be the ring of  germs of analytic functions on
$\matC^r$. We denote by $Mat_{(n,p)}({\cal{O}}_r)$ the set of all
matrices $n\times p$ with entries in ${\cal{O}}_r$. This set can be
identified with ${\cal{O}}^{np}_ r$, where ${\cal{O}}_r^{np}$ is a
free module of rank $np$.

We concentrate our attention in this paper to codimension $2$
determinantal singularities describing these singularities and their
deformations. The following proposition follows from the Auslander-
Buchsbaum  formula and the Hilbert- Burch's Theorem.

\begin{proposition}(\cite{FK}, pg. 3994)
\begin{itemize}
\item [1. ]Let $M$ be a $(n+1)\times n$ matrix  with entries on ${\cal{O}}_r$
and $f=(f_1,\ldots, f_{n+1})$ its maximal minors and, by abuse of
notation, the ideal generated by them. If $codim(V(f))\geq2$ the
following sequence
$$0\longrightarrow ({\cal{O}}_r)^n\longrightarrow
({\cal{O}}_r)^{n+1}\longrightarrow {\cal{O}}_r/(f)\longrightarrow
0$$ is exact. Moreover, ${\cal{O}}_r/(f)$ is Cohen-Macaulay and
$codim(V(f))=2$.
\item [2. ] If $X\subset\matC^r$ is Cohen-Macaulay, $codim(X)=2$ and
$X=V(I)$, then ${\cal{O}}_r/I$ has a minimal resolution of the type
$$0\longrightarrow {\cal{O}}_r^n\longrightarrow
({\cal{O}}_r)^{n+1}\longrightarrow {\cal{O}}_r/I\longrightarrow 0.$$
Moreover, there is an unity $u\in{\cal{O}}_r$ such that $I=u\cdot
f$, where $f$ is again the ideal of the  maximal minors of $M$.
\item[ 3.] Any deformation of $M$ is a deformation of
$X$;
\item[4. ] Any deformation of $X$ can be generated by a
perturbation of the matrix $M$.
\end{itemize}
\end{proposition}

It follows from this proposition that any deformation of a
Cohen-Macaulay  variety of codimension $2$ can be given as a
perturbation of the presentation matrix. Therefore, we can study
these varieties and their deformations using their representation
matrices. We can express the normal module and the space of the
first order deformations in terms of matrices, hence we can treat
the base of the semi-universal deformation using matrix
representation.


The singularity theory of $(n+1)\times n$ matrices has been studied
in \cite{FK} and \cite{Miriam}.

\section{Polar Multiplicity}

In this section we present the definition of polar varieties given
in \cite{tess}. This concept has been used by B. Tessier and L\^{e}
D. T. with the purpose of studying singularities of analytic
varieties.

Suppose that $f:X\longrightarrow S$ is a morphism of reduced complex
analytic spaces such that the fibers of $f$ are smooth of
codimension $d=dim X-dimS$, out of a closed, not dense set $F\subset
X$. Generally, we can embed $X\subset S\times \matC^N$ as in the
diagram

$$\xymatrix{ \ar[d]^f X \ar[r] & S\times\matC^N\ar[ld]^{p_1} \\
                S    &  }$$

\begin{definition}
Let $D_{d-k+1}$ be a subspace of $\matC^N$ of codimension $d-k+1$
and $G_d(\matC^N)$ the Grassmannian  of $d$-planes in $\matC^N$ . We
define
$$C(D_{d-k+1})=\{T\in G_d(\matC^N)\,|\,dim(T\cap D_{d-k+1})\geq k\}.$$
\end{definition}

Observe that since $T $ is a plane of dimension $d$, the expected
dimension of $T\cap D_{d-k+1}$ is $d-(d-k+1)=k-1$. Therefore, the
points in $C(D_{d-k+1})$ are points with exceptional contact with
$D_{d-k+1}$.

Let $p:\matC^N\longrightarrow\matC^{d-k+1}$ be a generic linear
projection whose kernel is $D_{d-k+1}$ (see \cite{tess}, pg. 419).
If $s\in X-F$, the fiber $X_s$ of $f$ in $X$ is non singular
contained in $\{f(s)\}\times\matC^N$. We denote by
$p_s:X_s\longrightarrow\matC^{d-k+1}$ the restriction to $X_s$ of
the projection $p$.

We define
$$P_k(f,p)^0=\{x\in X-F\,|\,x\in\Sigma(p_s)\}$$
and denote by $P_k(f,p)$ its closure in $X$.

\begin{proposition}(\cite{tess}, IV 1.3.2)
Suppose that $f:X\longrightarrow S$ is a flat map with smooth fibers
at every point of $X-F$. Then $P_k(f,p)$ is a closed analytic
subset, empty or of pure codimension $k$ em $X$.
\end{proposition}


\begin{definition}
Let $f:(X,0)\longrightarrow(S,0)$ be a morphism as above, a $S$-
embedding $(X,0)\subset(S,0)\times(\matC^N,0)$ and a generic linear
subspace $D_{d-k+1}\subset\matC^N$. The closed analytic subspace
$P_k(f,p)$ of $X$ is the relative polar variety of $X$ with
codimension $k$ associated to $f$ and $D_{d-k+1}$. Here $D_{d-k+1}$
denotes the kernel of $p$.
\end{definition}

When $S$ is a point, we call $P_k(f,p)$ \textit{absolute polar
variety}. In general, we omit the adjective ``relative" or
``absolute".

The key invariant of $P_k(f,D_{d-k+1})$ is its multiplicity,
$m_0(P_k(f,D_{d-k+1}))$, called the relative polar multiplicity of
$X$, and denoted by $m_k(X,f)$. If $f$ is the constant map, we
denote the polar multiplicity by $m_k(X)$. For a generic projection
the multiplicity is independent of $D_{d-k+1}$ and, in fact, is an
analytic invariant of $X$.

Let $X\subset\matC^N$ be a complex analytic variety of complex
dimension $d$. The absolute polar variety of  $X$ of codimension $d$
consists of a finite number of points or is empty. In both cases,
their multiplicity is not well defined. However, T. Gaffney in
\cite{Gaff1} introduces the $ d $-th polar variety by the following
construction:

Let $\mathfrak{X}\subset\matC^N\times\matC^s$ be  a complex analytic
variety of complex dimension $d+s$ and
$f:\mathfrak{X}\longrightarrow \matC^s$ an analytic function such
that $f^{-1}(0)=X$.  Then, we can define $m_d(X,p,f)=m_0(P_d(f,p))$,
where $P_d(f,p)$ is the polar variety of $\mathfrak{X}$ with respect
to $(f,p)$. In general, $m_d(X,p,f)$ depends on the choices of
$\mathfrak{X}$ and $p$, but when $\mathfrak{X}$ is a versal
deformation of $X$ or in the case that $X$ has a unique smoothing
(see the next section for the definition of the smoothing of an
analytic variety), $m_d$ depends only on $X$ and $p.$ Furthermore,
if $p$ is a generic linear embedding, $m_d$ is an invariant of the
analytic variety $X$, which we denote by $m_d(X).$

%

\section{The Generic Fiber}

Let $X_0\subset\matC^N$ be the germ of an analytic $d$-dimensional
variety, on some open set of $\matC^N$ with isolated singularity at
the origin. A smoothing of $X_0$ is a flat deformation with the
property that its generic fiber are smooth. More precisely:
\begin{definition}
We say that a germ of analytic variety $(X_0,0)$ with isolated
singularity of complex dimension $d\geq 1$ has a smoothing, if there
exist an open ball $B_{\epsilon}(0)\subset\matC^N$ centered at the
origin, a closed subspace $X\subset B_{\epsilon}(0)\times D$, where
$D\subset\matC$ is an open disc with center at zero and a proper
analytic map
$$F:X\longrightarrow D,$$
with the restriction to $X$ of the projection $p:B\times
D\longrightarrow D$ such that
\begin{itemize}
   \item[a)] $F$ is flat;
   \item[b)] $(F^{-1}(0),0)$ is isomorphic to $(X_0,0)$;
   \item[c)] $F^{-1}(t)$ is non singular for $t\neq0$.
\end{itemize}
\end{definition}

It follows from the above definition  that $X$ has isolated
singularity at the origin and is a normal variety if $X_0$ is normal
at zero. Moreover,
$$F\mid_{F^{-1}(D-0)}: F^{-1}(D-0)\longrightarrow D-\{0\}$$
is a fiber bundle whose fibers $X_t=F^{-1}(t)$ are non singular.

The topology of the generic fiber of a reduced curve has been
intensively studied (see \cite{Bg}). For instance, the following
result holds:

\begin{theorem}(\cite{Bg}, pg. 258)\label{confib}
Let $f:Y\longrightarrow D$ be a good representative of a flat family
$f: (Y, 0)\longrightarrow(D , O)$ of reduced curves. Then, for all
$t\in D$ the fiber $Y_t$ is connected.
\end{theorem}
\vspace{0.3cm}

For $n$- dimensional analytic spaces the following result is due to
Greuel and Steenbrink.

\begin{theorem}(\cite{greuel}, pg. 17)
Let $(X,0)$ be a complex analytic space, $n$-dimensional, with
isolated singularity and $X_t$ the Milnor fiber of a smoothing of
$(X,0)$. Then, $\Pi_i(X_t)=0$ for $i\leq dimX-codimX$.
\end{theorem}
\vspace{0.3cm}

It follows from the previous theorem that if $(X,0)$ is Cohen-
Macaulay of codimension $2$ with isolated singularity, then its
Milnor fiber is $(dimX-2)$- connected.

As a consequence of Sard's Theorem, it follows that complete
intersections are smoothable; moreover the base of their
semiuniversal deformations is smooth  whence the existence and
uniqueness of the smoothing hold for them. For determinantal
singularities, the existence and uniqueness of the smoothing do not
occur in general (see \cite{greuel}). But the following result was
proved by J. Wahl:

\begin{theorem}(\cite{Wahl}, pg. 241)
Let $(X,0)$ be a determinantal variety with isolated singularity at
the origin defined by $t\times t$ minors of a $n\times p$ matrix
$M$, whose entries are in ${\cal{O}}_r$, $2\leq t\leq n\leq p$. If
$dim(X)< n + p - 2t + 3$, then $X$ has a smoothing.
\end{theorem}
\vspace{0.3cm}

In particular, it follows from this result that if $(X,0)$ is Cohen-
Macaulay with codimension less than or equal to $2$ and
$dim(X,0)\leq 3$, then $(X,0)$ admits a smoothing. We also observe
that for Cohen- Macaulay singularities of codimension less than or
equal to $2$, there is no obstruction for lifting second-order
deformations, the basis of the semi-universal deformation is smooth
(\cite {FK}).

The following result was proved by Greuel and Steebring in
\cite{greuel}.

\begin{theorem}(See \cite{greuel}, pg. 540)\label{betti}
Let $X_t$ be the Milnor fiber of a smoothing of a normal
singularity, then $b_1(X_t)=0$.
\end{theorem}

\section{Morse Theory and the Topology of Varieties with Isolated Singularity}
Let $(X,0)\subset(\matC^N,0)$ be a $n$- dimensional variety with
isolated singularity at the origin. Suppose that $X$ has a
smoothing, i. e., there exist a flat family
$$\Pi:\mathfrak{X}\longrightarrow D\subset\matC,$$
restriction of the projection $\Phi:B_{\epsilon}(0)\times
D\longrightarrow D$, such that $X_t=\Pi^{-1}(t)$ is smooth for all
$t\neq0$ and $X_0=X$.

The variety $\mathfrak{X}$ also has isolated singularity at the
origin. Let $p$ be a complex analytic function defined in $X$ with
isolated singularity at the origin. Let
\begin{align*}
\widetilde{p}:\mathfrak{X}\subset\matC^N\times\matC&\longrightarrow\matC\\
(x,t)&\longrightarrow\widetilde{p}(x,t),
\end{align*}
such that $\widetilde{p}(x,0)=p(x)$ and for all $t\neq 0$, let
$\widetilde{p}(\,\cdot\,,t)=p_t$ be a Morse function in $X_t$.

Thus we have the following diagram

$$\xymatrix{&\hspace{.6cm} X_t\ar[d]^{p_t}\subset \hspace{.9cm}{\mathfrak{X}}  &\subset  \matC^N\times\matC\ar[d]^{(\Pi,p)} & \cr
            &\matC\times\{t\}   &  \matC\times\matC& }$$
Notice that the number of critical points of $p_t$ is finite. In
fact, $x$ is a critical point of $p_t$ if and only if $x$ is a
critical point of the function $Re(p_t):X_t\longrightarrow\matR$.
Since the real part of $p_t$ is an analytic function on $X_t$, the
number of critical points of $Re(p_t)$ and, hence of $p_t$, is
finite.

\begin{proposition}\label{saeki}
Let $X$ be a $n$- dimensional variety with isolated singularity at
the origin admiting a smoothing and $p_t:X_t\longrightarrow\matC$,
$p_t=\widetilde{p}(\,\cdot\,,t)$ as above. Then,
\begin{itemize}
\item[a)] If $t\neq0$
\begin{equation}\label{decomp}
X_t\simeq p_t^{-1}(0)\dot{\cup} \{\mbox{cells of dimension } n\},
\end{equation}
where $\dot{\cup}$ indicates the gluing of the spaces and $\simeq$
indicates that the spaces have the same homotopy type.
\item[b)] \begin{equation}\label{decomp1}
{\cal{X}}(X_t)={\cal{X}}((p_t)^{-1}(0))+(-1)^2n_{\sigma},
\end{equation}
where $n_{\sigma}$ is the number of critical points of $p_t$ and
${\cal{X}}(X_t)$ denotes the Euler characteristic of $X_t$.
\end{itemize}
\end{proposition}
\noindent \textbf{Proof.}

Let $x_1,...,x_{\nu}$ be the critical points of $p_t$ and
$y_i=p_t(x_i)$, $1\leq i\leq \nu$, their critical values. Suppose
that $0$ is a regular value of $p_t$, for all $t\neq0$. We denote by
$E_i$ the line segments connecting the points $y_i$ to $0$, $E_i\cap
E_j=\{0\}$ for $i\neq j$ and $E=\cup E_i$. Take $\epsilon>0$ small
enough such that $y_i\in D_{\eta}(0)$ for all $1\leq i\leq \nu$.

%
%
%
%
%
%
%
%
%
%
%
%
%
%
%
%
%
%
%
%
%
%
%
%
%
%
%

The set $D_{\eta}(0)$ is a regular neighborhood of $E$ that retracts
to $E$.
We can realize this retraction through a smooth vector field that
can be lifted into the stratified space $X_t$. Integrating this
vector field, the space $p_t^{-1}(D_{\eta})$ retracts by deformation
on $p_t^{-1}(E)$.

Then,
\begin{equation*}
X_t=p_t^{-1}(D_{\eta})\simeq p_t^{-1}(E)=\bigcup_i\left(
p_t^{-1}(E_i)\right)=p_t^{-1}(0)\cup \overline{\left(\bigcup_i
p_t^{-1}(E_i-\{0\})\right)}.
\end{equation*}

Observe first that $x_i $ is critical point of the restriction of
$p_t$ to $p_t^{-1}(E_i-\{0\})$ if and only if $x_i$ is a critical
point of the restriction of the real part of $p_t$ to
$p_t^{-1}({E_i-\{0\}})$. Therefore, it follows from the classical
Morse theory (see \cite{Milnor1}) that
\begin{equation}
X_t=p_t^{-1}(0)\cup \overline{\left(\bigcup_i
p_t^{-1}(E_i-\{0\})\right)}\simeq p_t^{-1}(0)\dot{\cup}
\{\mbox{cells of dimension } d\},
\end{equation}
where $\dot{\cup}$ indicates the gluing of the spaces and $\simeq$
indicates that the spaces have the same homotopy type.

As the Euler characteristic is a homotopy invariant, using the
decomposition (\ref{decomp}) we have
\begin{equation}\label{decomp1}
{\cal{X}}(X_t)={\cal{X}}((p_t)^{-1}(0))+(-1)^2n_{\sigma},
\end{equation}
where $n_{\sigma}$ is the number of critical points of $p_t$.
$\hfill \blacksquare$ \vspace{0.5cm}

A consequence of the decomposition $(\ref{decomp})$ is that only
$\overline{p_t^{-1}(E_i-\{0\})}$ contributes to the free part of
$H_n(X_t,\matZ)$. Hence, $b_n (X_t)$ is less than or equal to the
number of critical points of $p_t$.

\begin{remark}
This result also appears in  \cite{Kav}.
\end{remark}

\section{Determinantal Varieties}

In this section, we restrict our  attention to Cohen- Macaulay
singularities of codimension $2$ with isolated singularity at the
origin. These include determinantal surfaces in $\matC^4$ and $3$-
dimensional determinantal varieties in $\matC^5$. These varieties
admit a unique smoothing (see \cite{Miriam}), thus the following
definition makes sense:

\begin{definition}
Let $(X,0)\subset(\matC^r,0)$ be the germ of a codimension $2$
determinantal variety with isolated singularity at the origin,
$dim(X)=2,3$.The Milnor number of $X$, denoted $\mu(X)$, is defined
by $\mu(X)=b_d(X_t)$, where $X_t$ is the generic fiber of $X$ and
$b_d(X_t)$ is the $d$-th Betti number of $X_t$, $d=dim(X)$.
\end{definition}

Let $p:X\longrightarrow\matC$ be a complex analytic function with
isolated singularity at the origin. Then, $Y=X\cap p^{-1}(0)$ is a
variety of dimension $d-1$, with isolated singularity at $0$.

In particular, when $p:\matC^r\longrightarrow\matC$ is a linear
function, it follows from the presentation matrix that $Y$ is also a
Cohen- Macaulay determinantal variety of dimension $d-1$ on
$\matC^{r-1}$.

In the following theorem, we obtain a formula of type L\^{e}-Greuel
for germs of Cohen-Macaulay determinantal surfaces of codimension
$2$ with isolated singularity at the origin.


\begin{theorem}\label{polmil}
Let $(X,0)$ be the germ of a determinantal surface with isolated
singularity defined as above. Then,
$$m_2(X)=\mu(p^{-1}(0)\cap X)+\mu(X),$$
where $m_2(X)$ is the polar multiplicity of  $X$.
\end{theorem}
\noindent \textbf{Proof.}

From (\ref{decomp}), we have
\begin{equation}\label{decomp1}
{\cal{X}}(X_t)={\cal{X}}((p_t)^{-1}(0))+(-1)^2n_{\sigma},
\end{equation}
where $n_{\sigma}$ is the number of critical points of
$p_t:X_t\longrightarrow\matC$. It follows from de-finition of polar
multiplicity that $n_{\sigma}$ coincides with the $2$-th polar
multiplicity of the variety $P_k(f,p)$, i. e., $n_{\sigma}=m_2(X)$.

Moreover,
$${\cal{X}}(X_t)=b_0(X_t)-b_1(X_t)+b_2(X_t).$$
Then, using (\ref{decomp1}) we get
\begin{equation}
b_0(X_t)-b_1(X_t)+b_2(X_t)={\cal{X}}((p_t)^{-1}(0))+m_{2}(X).
\end{equation}


We know that $p_t^{-1}(0)\subset X_t$ is the generic fiber of the
determinantal curve $p_0^{-1}(0)\subset X$. Therefore, $p_t^{-1}(0)$
has the homotopy type of a bouquet of spheres of real dimension $1$.
Let $C=p_t^{-1}(0)\cap X_t$. Then, ${\cal{X}}(C)=1-\mu(C,0)$, where
$\mu(C,0)$ is the Milnor number of the curve (\cite{Bg}).

Since determinantal varieties are normal varieties, it follows from
Proposition \ref{betti} that $b_1(X_t)=0$. Moreover $X_t$ is
connected. Therefore,
$$1+b_2(X_t)=1-\mu(p_t^{-1}(0)\cap X_t)+m_2(X).$$
Hence
$$m_2(X)=\mu(p_t^{-1}(0)\cap X_t)+\mu(X).$$
$\hfill \blacksquare$ \vspace{0.8cm}

When $dim(X)=3$, we obtain an expression which reduces to
L\^{e}-Greuel formula when $b_2(X_t)=0$.

\begin{proposition}
Let $(X,0)\subset(\matC^5,0)$ be the germ of a determinantal variety
of codimension $2$ with isolated singularity at the origin. Then,
$$m_3(X)=\mu(p^{-1}(0)\cap X)+\mu(X)+b_2(X_t),$$
where $m_3(X)$ is the polar multiplicity of $X$.
\end{proposition}
\vspace{0.8cm}

\section{Index of $1$-Forms on Determinantal Varieties}

In this section we relate the formulas of the previous section with
Ebeling and Gusein-Zade index formulas in (\cite{EG}). They define
indices of $1$- forms on determinantal varieties having an essential
isolated singularity, EIDS. These singularities can be represented
by a matrix $M=(m_{ij}(x))$, $x\in\matC^r$, which is transverse,
away from the origin to the rank stratification of
$Mat_{(n,p)}(\matC)$, (see \cite{EG} for more details).

In particular, codimension two determinantal varieties with isolated
singularities are EIDS, and the results in \cite{EG} apply to this
class of singularities.

Let $X$ be a germ of codimension two determinantal variety with
isolated singularity and $\mathbf{\omega}$ the germ of a $1$-form on
$\matC^r$ whose the restriction to $(X,0)$ has an isolated singular
point at the origin.

Ebeling and Gusein-Zade definition of the Poincar\'e-Hopf index of
$\mathbf{\omega}$ reduces in our case to the following:

\begin{definition}
The Poincar\'e-Hopf index (PH-index), $ind_{PH}\omega$, is the sum
of the indices of the zeros of a generic pertubation
$\widetilde{\mathbf{\omega}}$ of the $1$-form $\mathbf{\omega}$ on
$\widetilde{X}$, a smoothing of $X$.
\end{definition}

\begin{proposition}(\cite{EG}, pg. 7)\label{PH-P}
The PH-index $ind_{PH}\omega$ of the $1$-form $\omega$ on the EIDS
$(X,0)$ is equal to the number of non-degenerate singular points of
a generic deformation $\widetilde{\omega}$ of the $1$-form $\omega$
on $\widetilde{X}_{reg}$, the regular part of $\widetilde{X}$.
\end{proposition}

For determinantal varieties with isolated singularity, the relation
between the PH-index and the radial index (see \cite{EG1} for the
definition of the radial index), is given by
$$ind_{PH}(\omega;X,0)=ind_{rad}(\omega;X,0)+(-1)^{dim(X)}\overline{{\cal{X}}}(X,0),$$
where $\overline{{\cal{X}}}(X,0)={\cal{X}}(X,0)-1$. The PH- index is
closely related to the $d$-th polar multiplicity as we can see in
the following proposition.
\begin{proposition}\label{ind}
Let $(X,0)\subset\matC^4$ be a determinantal surface, with isolated
singularity at the origin. Then, $ind_{PH}(\omega;X,0)=m_2(X)$,
where $\omega=dp$, here $p$ is a generic projection and $m_2(X)$ is
the second polar variety of $X$.
\end{proposition}
\noindent \textbf{Proof.}

Let $p:(X,0)\longrightarrow\matC$ be a generic linear projection and
$\epsilon>0$ small enough such that the restriction of $p$ to $X\cap
B_{\epsilon}(0)$ has isolated critical point at the origin. Then, by
Theorem $3$ of \cite{EG1}
$$ind_{rad}(dp;X,0)=(-1)^d\overline{{\cal{X}}}(p^{-1}(t)),$$
$t\neq 0$, where $\overline{{{\cal{X}}}(X)}={{\cal{X}}}(X)-1$. Then,
$$ind_{PH}(dp;X,0)=(-1)^2\overline{{\cal{X}}}(p^{-1}(t))+(-1)^{2}\overline{{\cal{X}}}(X,0).$$
Therefore,
$$ind_{PH}(dp;X,0)=({\cal{X}}(p^{-1}(t))-1)+({\cal{X}}(X,0)-1)=\mu(C)+\mu(X),$$
where $C=p^{-1}(0)\cap X$.

$\hfill \blacksquare$ \vspace{0.8cm}

\begin{remark}
  \begin{itemize}
     \item[a) ] This result is useful in calculations of
     $\mu(X)$, since in many cases one can use
     geometric methods to calculate $ind_{PH}(\omega;X,0)$. This
     procedure will be useful in the calculations in the next section.
     \item[b) ] An important problem not addressed in this
     work is the determination of an algebraic formula for the polar multiplicity as in
     L\^{e}-Greuel's formula for ICIS. See \cite{GaftNiv}, for
     an  algebraic approach characterizing the $d$- polar multiplicity
     of $d$- dimensional singular spaces.
  \end{itemize}
\end{remark}

\section{Examples}\label{6.3}

In this section, we compute the Milnor number $\mu(X)$ for some
normal forms of the simple determinantal surfaces $X$ in $\matC^4$
classified by   Fr\"uhbis-Kr\"uger, A. Neumer \cite{FKn}.

To calculate the Milnor number, we use the formula $m_2(X)
=\mu^{(1)}(X)+\mu(X)$ from Corollary \ref{polmil}, where $\mu^
{(1)}(X)$ is the Milnor number of the curve $X\cap p^{(-1)}(0)$,
where $p:X \rightarrow\matC$ is a generic linear projection.

Using corollary \ref{ind}, $m_2 (X) = ind_{PH}(\omega, X, 0) $.
Moreover, if $X$ is simple and $p$ is a generic linear projection,
$Y =X\cap p^{-1}(0) $ is a simple determinantal curve and its Milnor
number can be calculated or we can directly use the table of simple
curves in \cite{FK} pg. 4008-4009.

To find $ m_2 (X) $, or equivalently, $Ind_ {PH}(\omega, X, 0) $, we
can follow one of the following procedures:

\begin{itemize}
  \item [ a)] To use the algorithm proposed by W. Ebeling and S. M. Gusein-Zade in
  \cite{EG};
  \item [ b)] To calculate the number of non
  degenerate singular points of the linear form
  $\omega$ defined on a smoothing of $X$.
  \item [ c)] To obtain a perturbation $X_t$ of $X$ with singular points
  points $p_1,...,p_l$ and we use the fact that
  $$\mu(X)=\sum_{i=1}^{l}\mu(X_t,p_i).$$
\end{itemize}

We illustrate each one of these procedures in the examples below:

\begin{itemize}
\item {\underline{{\textbf{ $1^{\circ}$ Example }}}(see \cite{EG}, pg. 17):  Let $M=\left(
                         \begin{array}{ccc}
                           z & y & x \\
                           w & z & y \\
                         \end{array}
                       \right)$.
To apply Ebeling and Gusein-Zade  method let
$p:\matC^4\longrightarrow\matC$, $p(x,y,z,w)=w$ and
${\mathbf{\omega}}=dp$. We consider the space curve $(C,0)=X\cap
p^{-1}(0)$ represented by the matrix
$$N=\left(
    \begin{array}{ccc}
      z & y & x \\
      0 & x & y \\
    \end{array}
  \right).$$
The family
$$\left(
    \begin{array}{ccc}
      z & y+b & x+c \\
      a & x & y \\
    \end{array}
  \right),$$
is the versal unfolding of $(C,0)$, whose  discriminant is
$a(b^2-c^2)=0$ (see \cite{FK}).

We obtain $M$ from the versal deformation of $N$ taking $a=w$, $b=c
=0$ and, moreover, a smoothing $M_{\lambda}$ to $M$ is obtained
taking $a=b=w$, $c =\lambda \neq0$. For each fixed $\lambda$,
$M_{\lambda}$ intersects the discriminant in $3$ distinct points
where the function $P(x, y, z, w) = w$ has non-degenerate critical
points. Using \ref{polmil}, we obtain $\mu(X) +\mu(C)=3$ and
$\mu(C)=2$, it follows that $\mu(X)=1$.

\item\underline{\textbf{ $2^{\circ}$ Example:}}
Let $\left(
       \begin{array}{ccc}
         z & w+x & y^k \\
         w & y & x \\
       \end{array}
     \right).$
This normal form is equivalent to the second normal form in table
$2a$ in \cite{FKn}. Let $p:X \rightarrow\matC$ be defined by $p(x,
y, z, w)=w$, $\omega=dp$ and $(C, 0) $ the determinantal curve given
by
$$\left(\begin{array}{ccc}z &
x & y^{k} \\0 & y &x  \end{array}\right), k\geq1$$
whose versal unfolding is given by
$$\left(\begin{array}{ccc}z & x+b &
y^{k}+\displaystyle{\sum_{i=0}^{k-1}c_ix^i} \\a & y &x
\end{array}\right).$$

A smoothing of the determinantal surface $X$ is obtained by taking
$c_0=\lambda\neq 0$, $a=b=w$ and $c_i=0$ for $i\neq0$. Denote by
$M_{\lambda}$ the matrix obtained in this way.

For each $\lambda$ fixed, let
$f_{\lambda}:\matC^4\longrightarrow\matC^3$ be the map determined by
the maximal minors of $M_{\lambda}$, $Jf_{\lambda}$ the jacobian
matrix of $f_{\lambda}$, and $[Jf_{\lambda},\omega]$ the $4\times 4$
matrix whose first three rows are the rows of $Jf_ {\lambda}$ and
the last row are given by the coefficients of the form $\omega$.

To determine the number of non-degenerate critical points of
$\omega$ in the $M_{\lambda}$, we determine the number of solutions
of the equations

$$\begin{array}{l}
  zx-yw-w^2=0\\
  zy-wx^k-\lambda w=0\\
  y^2+yw-x^{k+1}-\lambda x=0\\
  z^2-kx^{k1}w^2=0\\
  zy+kx^kw=0\\
  wy+zx=0\\
  -kx^{k-1}w(2y+w)+z((k+1)x^k+\lambda)=0\\
  (k+1)x^ky+\lambda y=0\\
  2y^2+yw=0\\
  z(2y+w)-w((k+1)x^k+\lambda)=0\\
  (k+1)x^{k+1}+\lambda x=0\\
  2yx+xw=0.
\end{array}$$

The first tree equations above are determined by $2\times 2$ minors
of the matrix $M_{\lambda}$, the others are given by $3\times 3$
minors of $[Jf_{\lambda},\omega]$. In this case, we have $2k$
solutions of the form
$(x,(kx^{k+1})^{\frac{1}{2}},2kx^k,-2(kx^{k+1})^{\frac{1}{2}})$ and
$\lambda=(k+1)x^{k}$. We can verify that these solutions are non-
degenerate singular points.

Therefore, $ind_{PH}(\omega)=2k$. Then,
$$\mu(X)=2k-(k+1)=k.$$

\item\underline{\textbf{$3^{\circ}$ Example:}} Let $M=\left(
                                                       \begin{array}{ccc}
                                                         z & y & x \\
                                                         x & w & yz+y^kw \\
                                                       \end{array}
                                                     \right)
$. To determine the Milnor number for this normal form we use
procedure (c) and induction on $k$. We first consider the case
$k=1$. 

Let $p:\matC^4\longrightarrow\matC$ be given by $p(x,y,z,w)=y-z$ and
$\omega=dp$. In this case, the determinantal curve given by $C=X\cap
p^{-1}(0)$ is defined by the matrix
$$\left(
        \begin{array}{ccc}
          w & x & z^2 \\
          0 & z & x \\
        \end{array}
      \right)$$
that is the first normal form of \cite{FK} with $\mu(C)=3$. A
smoothing of $X$ is given by
$$\left(
        \begin{array}{ccc}
          w & x+\lambda & tz^2+yz+yw \\
          y & z & x \\
        \end{array}
      \right),$$
with $t\in\matC$, and $t\neq 0$.

We need to determine the number of common solutions of the equations
$$\begin{array}{l}
wz-y(x+\lambda)=0\\
wx-y(tz^2+yz+yw)=0\\
x(x+\lambda)-z(tz^2+yz+yw)=0\\
\end{array}$$
the defining equations of $X$, and of the equations determined by
the $3\times 3 $ minors of the matrix
$$[Jf_{\lambda},\omega]=\left(
    \begin{array}{cccc}
      -y             & -x-\lambda          & w          & z \\
       w             & -2yw-2yz-tz^2       & -y^2-2tzy  & x-y^2 \\
      2x-y^2+\lambda & -2y(x+\lambda)-z^2  & -2zy-3tz^2 & 0 \\
      0              & 1                   & -1         & 0 \\
    \end{array}
  \right).
$$

Solving this system of equations, we obtain the following solutions:
$$\begin{array}{l}
(0,0,0,\lambda)\\
(-\lambda,\pm\sqrt{-\lambda},0, w_i^+)\\
(-\lambda,\pm\sqrt{-\lambda},0, w_i^-),
\end{array}
$$
with $i=1,2$, $w_i^+$, $w_i^-$ are the solutions of equations
$$\begin{array}{l}
  w^2\mp2(-\lambda)w\pm\lambda\sqrt{-\lambda}=0\\
  w^2\pm2(-\lambda)w\mp\lambda\sqrt{-\lambda}=0.
\end{array}
$$

The other solutions satisfy $ x=\dfrac{y^2-\lambda}{2}$ and are
given by the equations
$$
\begin{array}{l}
  (1+3t)z^2+2yz+y(y^2+\lambda)=0\\
  w(\lambda+y^2)+2y(yz+tz^2)=0\\
  (y^2+\lambda)^2+4z(yz+tz^2)=0\\
 \end{array}
$$
It is possible to show in this case, that there are $3$ more
solutions of the form:
$$
  \left(\dfrac{y^2-\lambda}{2},y,\mu_3\sqrt[3]{\dfrac{(y^2+\lambda)(y^2+\lambda-2y)}{2(1+t)}},w(y)\right),
$$
where $w(y)$ is solution of the equation
$$\dfrac{y^2(y^2+\lambda)-y\mu_3^2\rho(y)(1+t)}{y^2+\lambda},$$
with
$$\rho(y)=\sqrt[3]{\dfrac{(y^2+\lambda)(y^2+\lambda-2y)}{2(1+t)}},$$
and $\mu_3$ is a $3$-rooth of the unit.

Then, $ind_{PH}(X;\omega,0)=8$. Therefore, $\mu(X)=5=\tau(X)-1$.

Suppose that for $k-1$, the Milnor number of $M$ is $2k+1$. To show
that for $k$, $\mu(X)=\tau(X)-1=2k+3$, we consider the following
$1$- parameter deformation of $M$
  $$M_t=\left(
      \begin{array}{ccc}
        z & y & x \\
        x & w & yz+y^kw+ty^{k-1}w \\
      \end{array}
    \right),
  $$
$t\in\matC$ e $t\neq 0$. The variety $X_t$ defined by the maximal
minors of $M_t$ is singular at the origin and on the points
$(0,-t,0,\pm\sqrt{-t^3})$. Then,
$$\mu(X)=\mu(X_t,0)+\mu(X_t,u_1)+\mu(X_t,u_2),$$
where $0$ is the origin in $\matC^4$,
$u_i=(0,-t,0,(-1)^i\sqrt{-t^3})$.

At the point $x=y=z=0$, $y\neq-t$, then $t+y$ is a unity on
${\cal{P}}$. Then,
$$
\begin{array}{llll}
&\left(\begin{array}{ccc}
      z & y & x \\
      x & w & yz+wy^{k-1}(t+y) \\
    \end{array} \right)&\sim
  &\left(\begin{array}{ccc}
      z & y & x/(t+y) \\
      x & w & \dfrac{yz}{t+y}+wy^{k-1} \\
    \end{array}\right)\sim \\
&\left(\begin{array}{ccc}
      z(t+y) & y & x/(t+y) \\
      x & w & yz+wy^{k-1} \\
    \end{array} \right)
&\sim
&\left(\begin{array}{ccc}
      z & y & x\\
      x & w & yz+wy^{k-1} \\
   \end{array}\right).
\end{array}
$$
By the induction hypotesis, $\mu(X_t,0)=2k+1$.

When $x=z=0$, $y=-t$ and $w=\pm\sqrt{-t^3}$, the matrix is
equivalent to
$$\left(\begin{array}{ccc}
      x & zy+ay^k+twy^{k-1} & w \\
      z & x & y \\
    \end{array} \right).$$
Calculating the $1$- jet of $zy+ay^k+twy^{k-1}$, we have
$$\left(\begin{array}{ccc}
      x & cy+dz & w \\
      z & x & y \\
    \end{array} \right)$$
with $c,\,d\neq 0$, that is therefore, equivalent to the first
normal form whose Milnor number is $1$.

Therefore, $\mu(X)=2k+1+1+1=2k+3$.

%
%
%
}
\end{itemize}

\begin{remark}
In these examples the equality $\tau(X)=\mu(X)+1$, holds where
$\tau$ is the Tjurina number. In \cite{Miriam} we verify that this
formula holds for almost all simple determinantal
surfaces in $\matC^4$ of Fr\"uhbis-Kr\"uger and Neumer list \cite{FKn}. We conjecture that
this formula holds more generally for all codimension $2$
Cohen-Macaulay isolated singularity. After completing this work, the
authors came across a paper of J. Damon and B. Pike \cite{JB} that
contains similar results. They discuss a more general inductive
procedure for computing the vanishing topology of matrix
singularities.
\end{remark}

$\mathbf{Acknowledgement.}$ We thank O. Saeki for suggesting us the
argument of Proposition \ref{saeki}.


\begin{thebibliography}{References}



\bibitem{FK} {\sc A. Fr\"uhbis-Kr\"uger}, {\em Classification of Simple Space Curves
Singularities,} {\sl Comm. in Alg.,} {\bf 27 (8)}, pp. 3993-4013,
(1999).

\bibitem{FKn} {\sc A. Fr\"uhbis-Kr\"uger, A. Neumer}, {\em Simple Cohen-Macaulay Codimension 2
Singularities,} {\bf arXiv:0808.2439v2}, (2008).


\bibitem{tess}{\sc B. Tessier}, {\em Variet$\acute{e}$s Polaires 2: Multiplicit$\acute{e}$s Polaires, Sections Planes, et Conditions de Whitney},
{\sl Actes de la conference de g$\acute{e}$ometrie
alg$\acute{e}$brique $\acute{a}$ la R$\acute{a}$bida }, {\sl
Springer Lecture Notes}, \textbf{ 961}, pp. 314- 491, (1981).

\bibitem{Gaff1} {\sc T. Gaffney}, {\em Polar Multiplicities and Equisingularity of Map Germs},
{\sl Topology}, \textbf{32}, pp. 185- 223, (1993).

\bibitem{Bg}
{\sc R.-O. G. Buchweitz and G.- M. Greuel }, {\it The Milnor Number
and Deformations of Complex Curve Singularities}, {\sl Inventiones
Mathematicae}, {\bf 58}, pp. 241- 281, (1980).










\bibitem{greuel} {\sc G. M. Greuel, J. Steenbrink}, {\it On the Topology of Smoothable Singularities}
{\sl Proceedings of Symposia in Pure Mathematics}, {\bf 40}, {\sl
Part 1}, pp. 535- 545, (1983).

\bibitem{GaftNiv} \sc{ T. Gaffney, N. Grulha Jr.}, {\it The multiplicity polar
theorem, collections of 1-forms and Chern numbers}, preprint.



%


\bibitem{JB} {\sc J. Damon and B. Pike}, {\it Solvable Groups, Free Divisors and Nonisolated Matrix Singulatities II: Vanishing Topology},
{\sl Available at http://www.math.unc.edu/\linebreak Faculty/jndamon/Nonisol.matrsing.II.v4.pdf}, (1981).

\bibitem{Wahl} {\sc J. Wahl}, {\it Smoothings of normal surface singularities},
{\sl Topology}, \textbf{20}, 219- 246, (1981).

\bibitem{Kav}{\sc K. Kaveh}, \it {Morse Theory and the Euler Characteristic of Sections of Spherical Varieties},
\sl{Transformation Groups}, {\textbf{9}}, {\bf{No. 1}}, pp. 47- 63,
(2004).




\bibitem{Bg}
{\sc R.-O. G. Buchweitz and G.- M. Greuel }, {\it The Milnor Number
and Deformations of Complex Curve Singularities}, {\sl Inventiones
Mathematicae}, {\bf 58},  {\sl 241- 281, (1980).}



\bibitem{Milnor1} {\sc W. J. Milnor}, {\it Morse Theory / Based on lecture notes by M. Spivak and R. Wells },
 {\sl Annals of Mathematics Studies}, \textbf{51}, {\sl New Jersey,
 (1963).}

\bibitem{Bruns} {\sc W. Bruns and U. Vetter}, {\it Determinantal
Rings}, {\sl Springer- Verlang}, \sl{New York}, (1998).


\bibitem{EG} \sc{ W. Ebeling and S. M. Gusein-Zade}, \it{On indices of $1$-forms on determinantal singularities},
\sl{Tr. Mat. Inst. Steklova}, \bf{267}, {\sl pp. 119- 131, (2009)}.

\bibitem{EG1} {\sc W. Ebeling and S. M. Gusein-Zade}, {\it Radial Index and Euler Obstruction of a $1$-form on a singular varieties},
{\sl Geometriae Dedicata}, {\bf 113}, {\sl pp. 231- 241, (2005).}

\bibitem{Miriam} {\sc M. S. Pereira}, {\it Variedades Determinantais e Singularidades de Matrizes},
{\sl Tese de Doutorado}, {\sl ICMC- USP, (2010).}

\bibitem{MC} {\sc M. A. S. Ruas and M. S. Pereira}, {\it The Milnor and Tjurina numbers of simple determinantal surface Singularities},
{\sl in preparation}.

\bibitem{BOT}\sc { J. J. NU\~{N}O-Ballesteros, B. Or\'efice, J. N. Tomazella}, {\it The Milnor of an Isolated Determinantal Variety}, preprint.



\end{thebibliography}
\end{document}